# Sécurité fonctionnelle : accorder la complexité des méthodes avec la complexité des systèmes

# Functional safety: matching the complexity of methods with the complexity of systems

F. Brissaud, L. O. Fernando, et B. Declerck
DNV France
69 rue Chevaleret, 75013 Paris, France
Tel. : +33 (0)1 83 94 05 01
E-mail : florent.brissaud@dnv.com

**Résumé**
En accord avec la norme de sécurité fonctionnelle CEI 61508, il est requis d'estimer l'intégrité de sécurité du système due aux défaillances aléatoires du matériel. Pour une fonction faiblement sollicitée, la mesure utilisée est la probabilité moyenne de défaillance dangereuse en cas de sollicitation ($PFD_{avg}$).
Dans la présente communication, quatre méthodes ont été appliquées à différentes configurations d'un cas d'étude : arbres de défaillance avec le logiciel GRIF/Tree, graphes de Markov multi-phase avec le logiciel GRIF/Markov, réseaux de Petri stochastiques à prédicats avec le logiciel GRIF/Petri, et équations approchées (développées par DNV et différentes de celles données dans la norme CEI 61508) avec le logiciel OrbitSIL.
Il est montré que toutes ces méthodes peuvent conduire à des résultats similaires pour l'estimation de la $PFD_{avg}$, en tenant compte des caractéristiques requises par la norme. Le choix d'une méthode doit alors se faire sans a priori, sur la base d'un accord entre les efforts de modélisation, les objectifs, et les propriétés du système. Afin d'assister l'analyste dans cette tâche, une discussion sur les avantages et limites de chacune de ces méthodes est présentée.

**Summary**
In line with the IEC 61508 functional safety standard, it is required to assess the safety integrity of a system due to random hardware failures. For a rarely used function (operating in a low demand mode), the measurement used is average probability of a dangerous failure on demand ($PFD_{avg}$).
In this paper, four methods have been applied to different configurations of a case study: failure tree analysis with the software GRIF/Tree, multi-phase Markov graphs with the software GRIF/Markov, stochastic Petri nets with predicates with the software GRIF/Petri, and approximate equations (developed by DNV and different from those given in the IEC 61508 standard) using the software OrbitSIL.
It is shown that all these methods can lead to similar results for the estimating of the $PFD_{avg}$, taking into account the required characteristics of the standard. The choice of method must be made without bias, based on an agreement between the modelling efforts, goals, and the system properties. To assist the analyst in this task, a discussion of the benefits and limitations of each of these methods is presented.

## I. Introduction

La norme de sécurité fonctionnelle CEI 61508 (CEI, 2010) (dont la seconde édition est disponible depuis 2010) propose une approche générique de toutes les activités liées au cycle de vie de sécurité de systèmes électriques et/ou électroniques et/ou électroniques programmables (E/E/PE) qui sont utilisés pour réaliser des fonctions de sécurité. Des normes d'application sectorielle sont basées sur cette dernière, et notamment la CEI 61511(CEI, 2004) (actuellement en cours de révision) pour les systèmes instrumentés de sécurité (SIS) dans le secteur des industries de transformation.

Une fois que les exigences de sécurité ont été allouées et spécifiées, en termes de fonctions de sécurité (SIF) et d'intégrité de sécurité (probabilité moyenne d'exécution satisfaisante des SIF), il est requis (entre autres) d'estimer l'intégrité de sécurité due aux défaillances aléatoires du matériel. Pour une SIF faiblement sollicitée (réalisée uniquement sur sollicitation, et pas plus d'une fois par an), la mesure utilisée est la probabilité moyenne de défaillance dangereuse en cas de sollicitation ($PFD_{avg}$), calculée par l'indisponibilité moyenne de sécurité. Il est alors requis que cette probabilité soit inférieure à un objectif chiffré de défaillance, défini par les spécifications et le niveau d'intégrité de sécurité (SIL) alloué.

Il est prescrit dans la norme CEI 61508 que le calcul de la $PFD_{avg}$ doit tenir compte de certaines caractéristiques : architecture du SIS et de ses sous-systèmes (e.g. « *M*-sur-*N* ») ; taux de défaillance de chaque sous-système et élément (pour les défaillances dangereuses détectées et non-détectées) ; défaillances de cause commune (e.g. avec le modèle du facteur *β*) ; temps de réparation ; intervalles de temps des essais périodiques (i.e. tests de révision) permettant de révéler les défaillances non-détectées, ainsi que l'efficacité de ces tests. Pour cette évaluation, plusieurs méthodes sont proposées, incluant les arbres de défaillance (Dutuit *et al.*, 2008 ; Signoret *et al.*, 2007) (mathématiquement équivalent aux blocs diagrammes de fiabilité), les graphes de Markov (multi-phase) (Bukowski 2005 & 2006, Signoret *et al.*, 2007 ; Dutuit *et al.*, 2008), les réseaux de Petri (stochastiques à prédicats) (Dutuit *et al.*, 2008 ; Signoret *et al.*, 2007), ainsi des équations approchées (Oliveira, 2008 ; Oliveira & Abramovitch 2010 ; Vaurio 2011 ; Brissaud *et al.* 2010).

Enfin, il est précisé qu'il appartient à l'analyste de déterminer laquelle de ces méthodes est la plus appropriée compte tenu des circonstances. Cette communication a pour vocation de fournir des éléments permettant d'assister ce choix.

Pour cela, un cas d'étude est tout d'abord présenté dans la Section II, considérant différents jeux de données. Dans la Section III, les arbres de défaillance, modèles Markoviens, réseaux de Petri, et équations approchées sont respectivement appliqués aux différentes configurations de ce cas d'étude. Les résultats sont ensuite présentés dans la Section IV, avec une synthèse comparative des méthodes. La Section V présente les conclusions qui closent cette communication.

# II. Cas d'Étude

## II.1. Description du système et des hypothèses

L'architecture du système est définie par une configuration en « *M*-sur-*N* » (en anglais « *M-out-of-N* »), c'est-à-dire que le système est composé de *N* « canaux » (i.e. sous-systèmes) et est capable de réaliser sa fonction si *M* canaux ou plus (n'importe lesquels parmi les *N*) sont dans un état opérant. Ce système est en fait une sous-partie (i.e. capteurs, unité de traitement, ou éléments finaux) d'un système instrumenté de sécurité (SIS) et, afin d'évaluer la SIF dans sa globalité, il doit être complété par les autres sous-parties.

La fonction du système est une fonction de sécurité qui est réalisée uniquement sur sollicitation, et dont la fréquence des demandes n'est pas plus grande qu'une fois par an. Cette fréquence de « demandes non-désirées » est assez faible pour ignorer les effets de ces demandes sur la disponibilité du système (en pratique, cette fréquence est communément comprise entre 1 tous les 10 ans et 1 tous les 1000 ans).

Chaque canal qui compose le système est considéré comme étant opérant si et seulement s'il n'est pas dans un mode de défaillance dangereux. Pour chaque canal, trois modes de défaillance (dangereux) sont considérés :
- défaillance dangereuse détectée en ligne (i.e. dès qu'elle se produit) par des auto-diagnostics ;
- défaillance dangereuse uniquement détectée par des tests de révision (périodiques), (et également détectée par des demandes réelles) ;
- défaillance dangereuse uniquement détectée par des demandes réelles (car les tests de révision ne sont pas efficaces à *100*%).

Chaque mode de défaillance se produit selon un taux de défaillance constant (i.e. distribution exponentielle), et les taux de défaillance sont les mêmes pour tous les canaux (mais, bien sûr, dépendent des modes de défaillance). Il est considéré que les ressources de maintenance requises sont toujours disponibles, c'est-à-dire qu'aucune logistique de délais de maintenance et qu'aucune indisponibilité d'équipe de maintenance n'est prise en compte. Les effets des durées de tests sont considérés comme négligeables pour l'intégrité de sécurité (et ne sont donc pas quantifiés).

L'intégrité de sécurité du système est calculée sur la base d'une période de temps qui est, par exemple, l'intervalle entre deux tests complets (qui ne peuvent être réalisés que par une demande réelle) ou, à défaut, la durée de vie planifiée du système.

## II.2. Notations

*MooN* architecture du système i.e. le système est composé de *N* « canaux » et est capable de réaliser sa fonction si *M* canaux ou plus sont dans un état opérant

$\lambda_D$  taux de défaillance de n'importe quel canal (qui compose le système), au regard des défaillances dangereuses

$\lambda_{DD}$  taux de défaillance de n'importe quel canal, au regard des défaillances dangereuses détectées en ligne par des auto-diagnostics

$\lambda_{DU}$  taux de défaillance de n'importe quel canal, au regard des défaillances dangereuses non détectées en ligne

$\lambda_{DUT}$  taux de défaillance de n'importe quel canal, au regard des défaillances dangereuses uniquement détectées par des tests de révision (et également détectées par des demandes réelles)

$\lambda_{DUU}$  taux de défaillance de n'importe quel canal, au regard des défaillances dangereuses uniquement détectées par des demandes réelles

DC  couverture des auto-diagnostics, telle que $\lambda_{DD} = DC \times \lambda_D$ et $\lambda_{DU} = (1 - DC) \times \lambda_D$

PTC  couverture des tests de révision, telle que $\lambda_{DUT} = PTC \times \lambda_{DU}$ et $\lambda_{DUU} = (1 - PTC) \times \lambda_{DU}$

$\beta_{DD}$  facteur de cause commune de défaillances, au regard des défaillances dangereuses détectées en ligne par des auto-diagnostics

$\beta_{DUT}$  facteur de cause commune de défaillances, au regard des défaillances dangereuses uniquement détectées par des tests de révision

$\beta_{DUU}$  facteur de cause commune de défaillances, au regard des défaillances dangereuses uniquement détectées par des demandes réelles

$\mu_{DD}$  taux de réparation de n'importe quel canal, au regard des défaillances dangereuses détectées en ligne par des auto-diagnostics

$\mu_{DUT}$  taux de réparation de n'importe quel canal, au regard des défaillances dangereuses uniquement détectées par des tests de révision

$T_1$  période de tests de révision
$T_0$  période de calcul pour l'intégrité de sécurité du système

$PFD_{avg}$  probabilité moyenne de défaillance dangereuse en cas de sollicitation (calculée sur la période $T_0$)

Chaque taux de défaillance peut alors être divisé en deux parties : défaillances indépendantes (non dues à des causes communes), et défaillances dues à des causes communes. Ces parties sont respectivement :
- $(1 - \beta_{DD}) \times \lambda_{DD}$ et $\beta_{DD} \times \lambda_{DD}$ pour les défaillances dangereuses détectées en ligne par des auto-diagnostics ;
- $(1 - \beta_{DUT}) \times \lambda_{DUT}$ et $\beta_{DUT} \times \lambda_{DUT}$ pour les défaillances dangereuses uniquement détectées par des tests de révision ;
- $(1 - \beta_{DUU}) \times \lambda_{DUU}$ et $\beta_{DUU} \times \lambda_{DUU}$ pour les défaillances dangereuses uniquement détectées par des demandes réelles.

**II.3. Jeux de données**

Afin de comparer différentes configurations d'un cas d'étude, les six jeux de données présentés dans le Tableau 1 sont considérés.

**Table 1.** Jeux de données pour le cas d'étude

|  | cas *i* | cas *ii* | cas *iii* | cas *iv* | cas *v* | cas *vi* |
|---|---|---|---|---|---|---|
| ***MooN*** | 1oo1 | | 1oo2 | | 2oo3 | |
| $\lambda_D$ [heure$^{-1}$] | $2.70 \times 10^{-6}$ | $1.35 \times 10^{-5}$ | $2.70 \times 10^{-6}$ | $1.35 \times 10^{-5}$ | $2.70 \times 10^{-6}$ | $1.35 \times 10^{-5}$ |
| ***DC*** | 0.50 | 0.25 | 0.50 | 0.25 | 0.50 | 0.25 |
| ***PTC*** | 0.90 | 0.70 | 0.90 | 0.70 | 0.90 | 0.70 |
| $\beta_{DD}$ | 0.02 | 0.05 | 0.02 | 0.05 | 0.02 | 0.05 |
| $\beta_{DUT}$ | 0.05 | 0.10 | 0.05 | 0.10 | 0.05 | 0.10 |
| $\beta_{DUU}$ | 0.05 | 0.10 | 0.05 | 0.10 | 0.05 | 0.10 |
| $\mu_{DD}$ [heure$^{-1}$] | 0.0417 | 0.0833 | 0.0417 | 0.0833 | 0.0417 | 0.0833 |
| $\mu_{DUT}$ [heure$^{-1}$] | 0.0417 | 0.0833 | 0.0417 | 0.0833 | 0.0417 | 0.0833 |
| $T_1$ [heures] | 4,383 | 8,766 | 4,383 | 8,766 | 4,383 | 8,766 |
| $T_0$ [heures] | 70,128 (i.e. 8 ans) | | | | | |

# III. Application des méthodes

**III.1. Arbres de défaillance avec GRIF/Tree (ALBISIA)**

Les arbres de défaillance consistent en une approche déductive qui exprime un évènement sommet par des combinaisons d'évènements de base, via des portes logiques telles que « ou », « et », et « *k*-sur-*n* ». Les analyses sont ensuite effectuées en utilisant de l'algèbre Booléenne, communément par l'intermédiaire des coupes minimales (MCS pour, en anglais, « *minimal cut sets* »), (i.e. ensemble minimal d'évènements de base dont l'occurrence/présence assure l'occurrence de l'évènement sommet). Parmi les cautions à rappeler pour les analyses par arbres d'évènement, il convient de noter que : la combinaison de probabilités moyennes d'évènements de base n'est pas équivalente à la probabilité moyenne de l'évènement sommet, (la probabilité de l'évènement sommet doit d'abord être exprimée en fonction du temps, et c'est ensuite que la moyenne peut être calculée) ; et si des troncatures de coupes minimales sont utilisées (afin de gagner du temps d'analyse pour les gros systèmes), l'ensemble des coupes minimales qu'il convient de « garder » peut être fonction du temps (Dutuit *et al.*, 2008). Pour effectuer des analyses par arbres de défaillance, le choix de l'outil logiciel est donc important. Dans la présente communication, le module « Tree » de la suite GRIF (GRIF, 2012) est utilisé. Ce logiciel est basé sur ALBIZIA, un calculateur développé par Total, qui utilise les diagrammes de décision binaire (en anglais, « *Binary Decision Diagram* » (BDD)).

Une caractéristique intrinsèque des arbres de défaillance est que tous les évènements de base sont indépendants (exceptés, bien sûr, les évènements de base répétés). Cela signifie notamment que si un évènement de base se produit, il ne peut pas prévenir l'occurrence d'un autre évènement de base. Par exemple, un mode de défaillance d'un canal ne prévent pas les autres modes de défaillance du même canal (i.e. des états multiples sont possibles au même instant), et une défaillance indépendante ne prévent pas une défaillance de cause commune (et réciproquement). Ceci peut conduire à une surestimation de la probabilité de l'évènement sommet si les évènements de base sont considérés comme incompatibles, et notamment s'ils sont nombreux (e.g. nombreux modes de défaillance). De plus, les arbres de défaillance (classiques) sont des modèles « statiques », c'est-à-dire que l'architecture du modèle (e.g. les portes logiques) ne dépendent pas du temps et/ou de variables aléatoires. Cependant, afin de modéliser des évènements sommets qui dépendent de l'ordre d'occurrence (i.e. la séquence) des évènements de base, des arbres de défaillance « dynamiques » ont aussi été développés (Cepin & Mavko, 2002).

L'arbre de défaillance appliqué au cas d'étude avec une architecture *1oo2* (cas *iii* et *iv*, cf. Table 1) est donné sur la Figure 1. L'évènement sommet représente la défaillance de la fonction de sécurité au niveau du système. Ensuite, plusieurs portes représentent les défaillances possibles des canaux qui composent le système. Au niveau le plus bas, se trouve les évènements de base qui représentent les défaillances indépendantes (un par canal) et de causes communes (un pour tous les canaux), au regard des défaillances dangereuses détectés en ligne par des auto-diagnostics, uniquement détectées par des tests de révision, et uniquement détectées par des demandes réelles (soit un total de neuf évènements de base).

En utilisant GRIF/Tree, une porte avec un triangle sur la droite signifie un transfert vers la porte correspondante avec un triangle en bas (e.g. Portes C2_D, sur la Figure 1). De plus, un évènement de base coloré en gris clair est une répétition de l'évènement de base correspondant coloré en gris foncé (e.g. Évènement de base 2, 4, et 6, sur la Figure 1). Un évènement de base et sa répétition sont physiquement les mêmes évènements. Enfin, la porte sommet de sur la Figure 1 est une porte « *M-sur-N* ». Sur cette figure, « 2/2 » signifie « *2-sur-2* » (ce qui est équivalent à une porte « et »). Comme les arbres de défaillance sont orientés « défaillances » et non « succès », cette porte est appropriée pour modéliser une architecture *1oo2*. Pour les évènements de base, trois types de lois sont ici utilisés :
- « GLM », qui correspond à des défaillances qui sont détectées en ligne, définie par une probabilité de défaillance au démarrage (non utilisée ici et donc égale à *0*), un taux de défaillance, et un taux de réparation ;
- « periodic-test », qui correspond à des défaillances qui sont périodiquement détectées, définie par un taux de défaillance, un taux de réparation, l'intervalle de temps entre deux tests consécutifs, et la date du premier test (ici égale au précédent paramètre, autrement des tests asynchrones peuvent être considérés) ;
- « exponential », qui correspond à des défaillances qui ne sont jamais détectées, définie par un taux de défaillance.

D'autres lois sont également disponibles dans GRIF/Tree, qui permettent notamment de modéliser d'avantages de caractéristiques telles que des lois de Weibull, des tests et/ou des réparations imparfaites, des défaillances générées par des tests et/ou des réparations, des indisponibilités et/ou des modes dégradés durant des tests, etc.

À noter que les taux de réparation des défaillances de causes communes sont égaux au nombre de canaux multiplié par le taux de réparation d'un canal (il est considéré que les ressources de maintenance sont toujours disponibles, cf. Section II.1). En effet, comme tous les évènements de base sont indépendants, lorsqu'au moins un canal est réparé, l'évènement de base qui représente la défaillance commune de tous les canaux n'est plus « valide ».

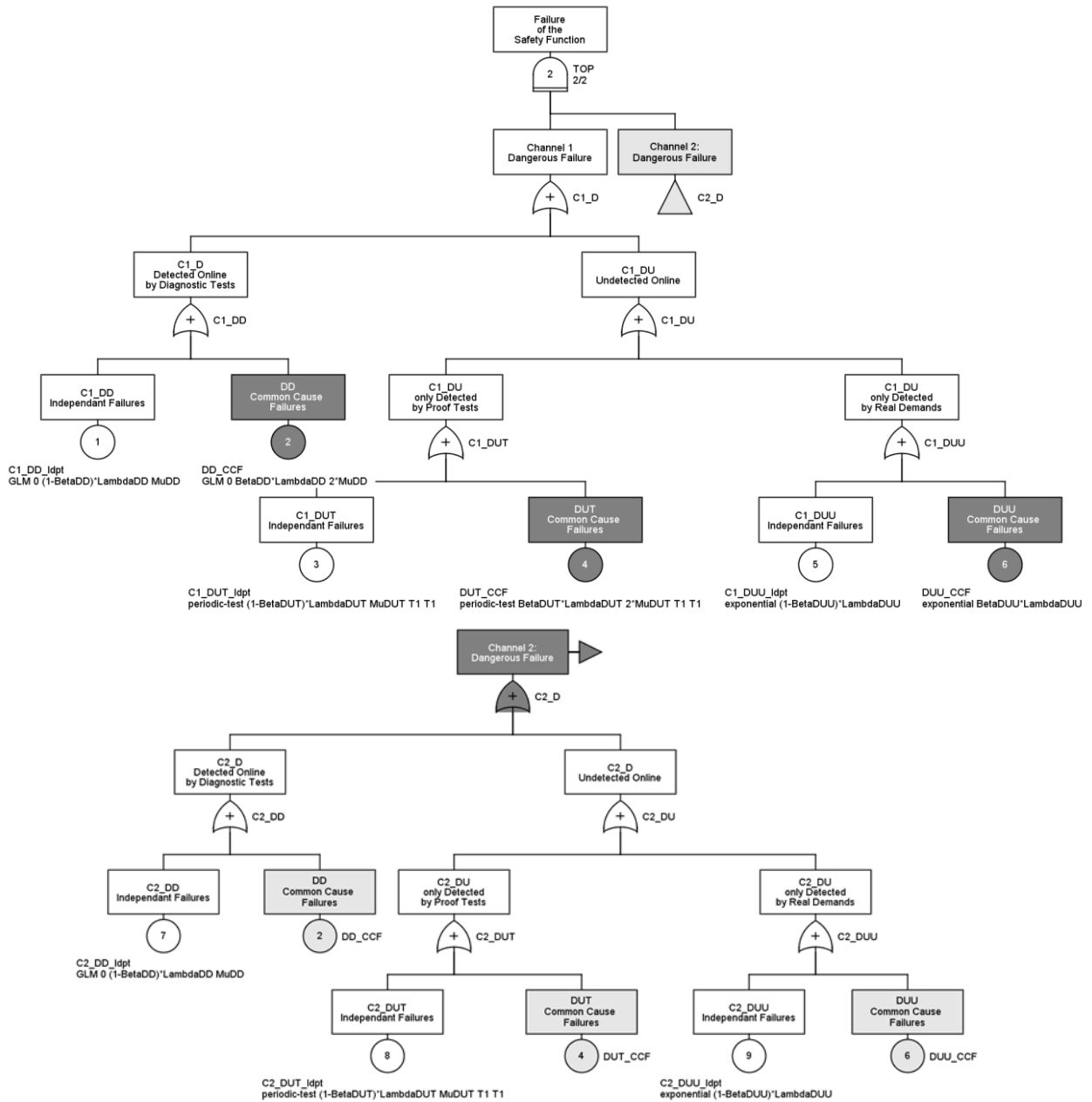

**Figure 1.** Arbres de défaillance avec GRIF/Tree, cas d'étude avec une architecture *1oo2*

### III.2. Graphes de Markov multi-phase avec GRIF/Markov (MARK-XPR)

Les modèles Markoviens sont des approches basés sur des états, qui permettent d'effectuer des analyses mathématiques comme, par exemple, le calcul des probabilités d'être dans des états donnés en fonction du temps, les fréquences moyennes d'entrées et de sorties dans des états donnés, et les temps moyens de séjours.

Un graphe de Markov est constitué d'états (cercles) et de transitions (flèches). Les états sont définis au niveau du système. Pour les systèmes à plusieurs composants, les états sont donc définis par les combinaisons possibles d'états des composants du système. Le nombre d'états dans un graphe de Markov est donc dépendant du nombre de composants et du nombre de modes (e.g. opérants, de défaillance, de réparation) à considérer pour chaque composant. Pour réduire le problème d'explosion combinatoire, il convient souvent de considérer que chaque composant ne peut pas se trouver dans de multiples états au même instant (e.g. un mode de défaillance prévient l'occurrence des autres modes de défaillance du même composant, contrairement aux arbres de défaillance, cf. Section III.1). De plus, il est souvent aussi souhaitable de regrouper plusieurs combinaisons dans un même état. Par exemple, si tous les composants du système sont identiques (e.g. mêmes taux de défaillance et taux de réparation), il peut être possible de raisonner en termes de nombre total de composants dans chaque mode, à la place de raisonner en termes de mode spécifique pour chaque composant. Au temps $t_0 = 0$, l'état initial est défini par des probabilités (e.g. l'état où tous les composants du système sont opérants est souvent défini comme l'état initial avec une probabilité de *1*). Les transitions entre états (e.g. défaillances, réparations) sont ensuite modélisées par des taux de transition, qui doivent être constants (i.e. distributions exponentielles). L'hypothèse Markovienne qui résulte de cela est que l'état au temps $t + \Delta t$ dépend uniquement de l'état au temps *t*, mais pas du temps *t*, ni des états précédents le temps *t*.

Parce que les tests périodiques se produisent à des instants déterministes (et donc ne suivent pas des distributions exponentielles), des modèles de Markov multi-phase ont été développées. Le système est alors modélisé au travers ses différentes phases (e.g. opération, test, réparation, etc.), en utilisant un graphe pour chaque phase. La séquence des phases est ensuite déterminée (cycliquement), tout comme la durée (constante) de chaque phase. De plus, lors du passage d'une phase à la suivante, une matrice (probabiliste) de liaison est utilisée pour définir le prochain état initial (dans la prochaine phase) en fonction de l'état courant (dans la phase courante). Dans la présente communication, le module « Markov » de la suite GRIF (GRIF, 2012) est utilisé. Ce logiciel est basé sur MARK-XPR, un calculateur développé par Total.

Le graphe de Markov multi-phase appliqué au cas d'étude avec une architecture *1oo2* (cas *iii* et *iv*, cf. Table 1) est donné sur la Figure 2. L'État 1 est l'état initial (représenté en couleur), où les deux canaux du système sont dans des modes opérants (nommés « OK »). Les autres états sont définis par les combinaisons possibles de modes de défaillance (nommés « DD », « DUT », et « DUU », suivant les notations données dans la Section II.2) et de modes de réparation (nommés « RepDUT » pour les modes de réparation « DUT »), (les modes de réparation « DD » sont les mêmes modes que les modes de défaillance « DD » car ces défaillances sont détectées en ligne, et les modes de réparation « DUT » n'ont pas besoin d'être modélisés d'après les hypothèses de la Section II.1). Puisque les canaux sont identiques (i.e. mêmes modes et mêmes transitions), les états du système sont définis indépendamment du canal qui cause le mode opérant, de défaillance, ou de réparation (e.g. « OK_DD » correspond au cas où un canal est dans un mode opérant et un canal dans un mode de défaillance « DD », indépendamment de quel canal est dans lequel de ces modes). En utilisant GRIF/Markov, le paramètre « Eff. » qui figure sous chaque état est la propriété à calculer. Dans ce cas, il s'agit de la probabilité de défaillance dangereuse en cas de sollicitation, qui est égale à *0* si au moins un canal est dans un mode opérant, et à *1* sinon (en accord avec l'architecture *1oo2*).

Pour ce cas d'étude, seule une phase est requise. Sa durée est égale à la période de tests de révision. Une fois que cette période est achevée, la phase redémarre. De plus, à chaque « redémarrage », la matrice de liaison est utilisée pour transférer chaque mode de défaillance détectée par un test de révision (i.e. mode de défaillance « DUT »), dans le mode de réparation correspondant (i.e. mode de défaillance « RepDUT »). Ici, cette matrice de liaison utilise uniquement des probabilités égales à *0* ou à *1*. Cependant, l'utilisation de matrices avec des probabilités comprises entre *0* et *1* et/ou l'utilisation de plusieurs phases permettent de modéliser d'avantages de caractéristiques telles que des tests asynchrones, des tests et/ou des réparations imparfaites, des défaillances générées par des tests et/ou des réparations, des indisponibilités et/ou des modes dégradés durant des tests, etc. La taille du modèle peut, néanmoins, augmenter drastiquement lorsque des caractéristiques supplémentaires sont prises en compte, jusqu'à rendre le modèle inexploitable dans certains cas.

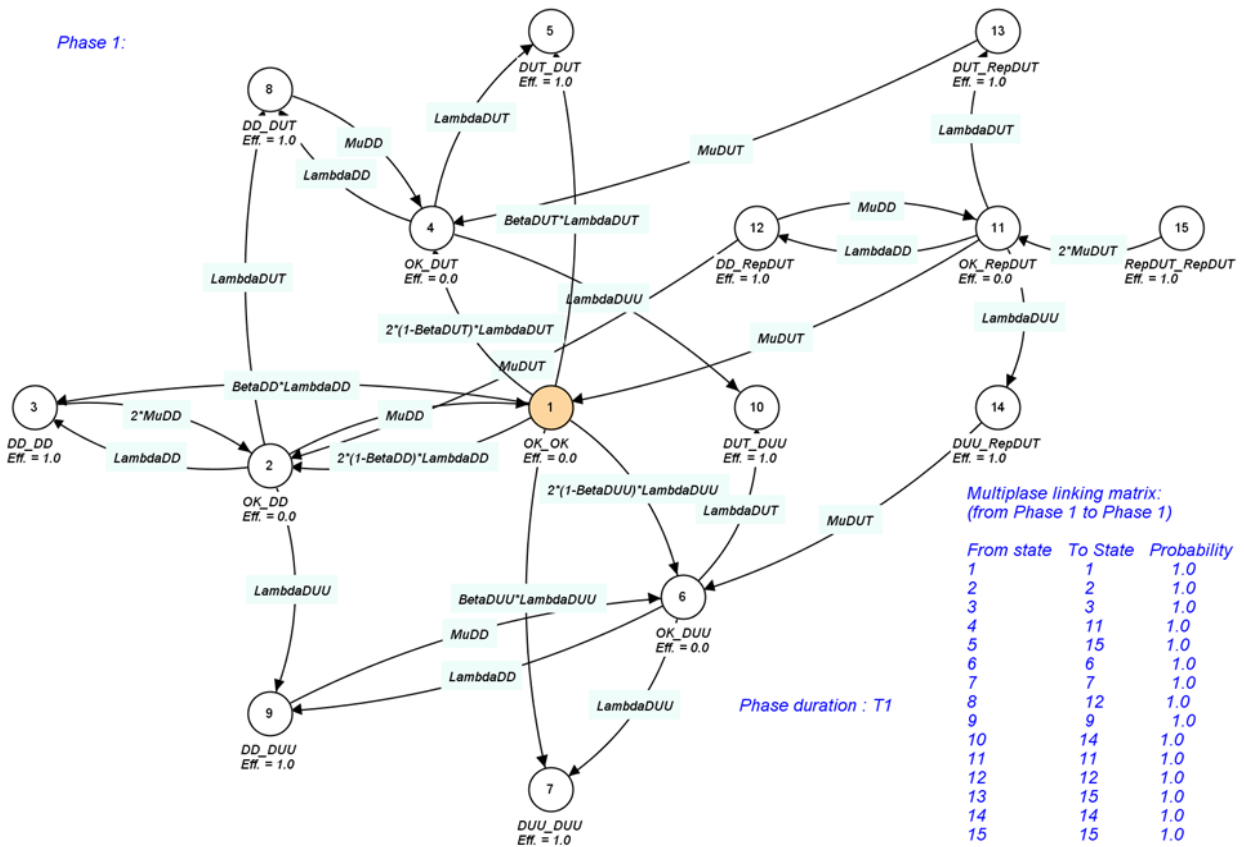

**Figure 2.** Graphes de Markov multi-phase avec GRIF/Markov, cas d'étude avec une architecture *1oo2*

### III.3. Réseaux de Petri stochastiques à prédicats avec GRIF/Petri (MOCA-RP V13)

Les réseaux de Petri fournissent un outil graphique pour modéliser les comportements (dynamiques) des systèmes, et ensuite permettre d'effectuer des analyses de disponibilité par simulations de Monte Carlo (i.e. les résultats sont obtenus statistiquement à l'issue de plusieurs histoires simulées).

Un réseau de Petri (classique) est constitué de places (cercles) et de transitions (rectangles). Des connections (arcs orientés) peuvent relier une place à une transition (arc d'entrée) ou vice-versa (arc de sortie), et peuvent être « valuées » (autrement, la valeur est de *1* par défaut). Les places peuvent contenir des jetons (petits cercles pleins) qui sont « déplacés » par l'intermédiaire des transitions lorsque celles-ci sont franchies. Une transition est franchissable lorsque chacune de ses places d'entrée (qui sont liées à la transition par un arc d'entrée) contient un nombre de jetons égal ou supérieur à la valeur de l'arc correspondant. Franchir une transition se fait alors en deux étapes : premièrement, le retrait dans chaque place d'entrée d'un nombre de jetons égal à la valeur de l'arc d'entrée correspondant ; deuxièmement, le dépôt dans chaque place de sortie d'un nombre de jetons égal à la valeur de l'arc de de sortie correspondant.

Généralement, les places d'un réseau de Petri représentent les objets ou conditions, les jetons précisent les valeurs de ces objets ou conditions, et les transitions modélisent l'activité du système. La dimension temporelle est introduite par des délais de franchissement de transitions (pendant lequel les transitions doivent aussi rester franchissables). Dans les réseaux de Petri stochastiques, ces délais sont des variables aléatoires. De plus, les réseaux de Petri à prédicats utilisent des variables pour inclure deux autres propriétés : des « gardes », variables ou expressions Booléennes qui rendent infranchissables les transitions tant qu'elles ne sont pas vérifiées ; et des « affectations », attributions qui modifient les valeurs de variables lors du franchissement des transitions. Dans la présente communication, le module « Petri » de la suite GRIF (GRIF, 2012) est utilisé. Ce logiciel est basé sur MOCA-RP, un simulateur de Monte Carlo développé par Total.

Le réseau de Petri appliqué au cas d'étude avec une architecture *1oo2* (cas *iii* et *iv*, cf. Table 1) est donné sur la Figure 3 (avec les jetons tels que définis au temps $t_0 = 0$). Les Places 1 à 5 modélisent le premier canal du système, les Places 12 à 16 modélisent le second canal du système, et les Places 6 à 11 modélisent les défaillances de causes communes.

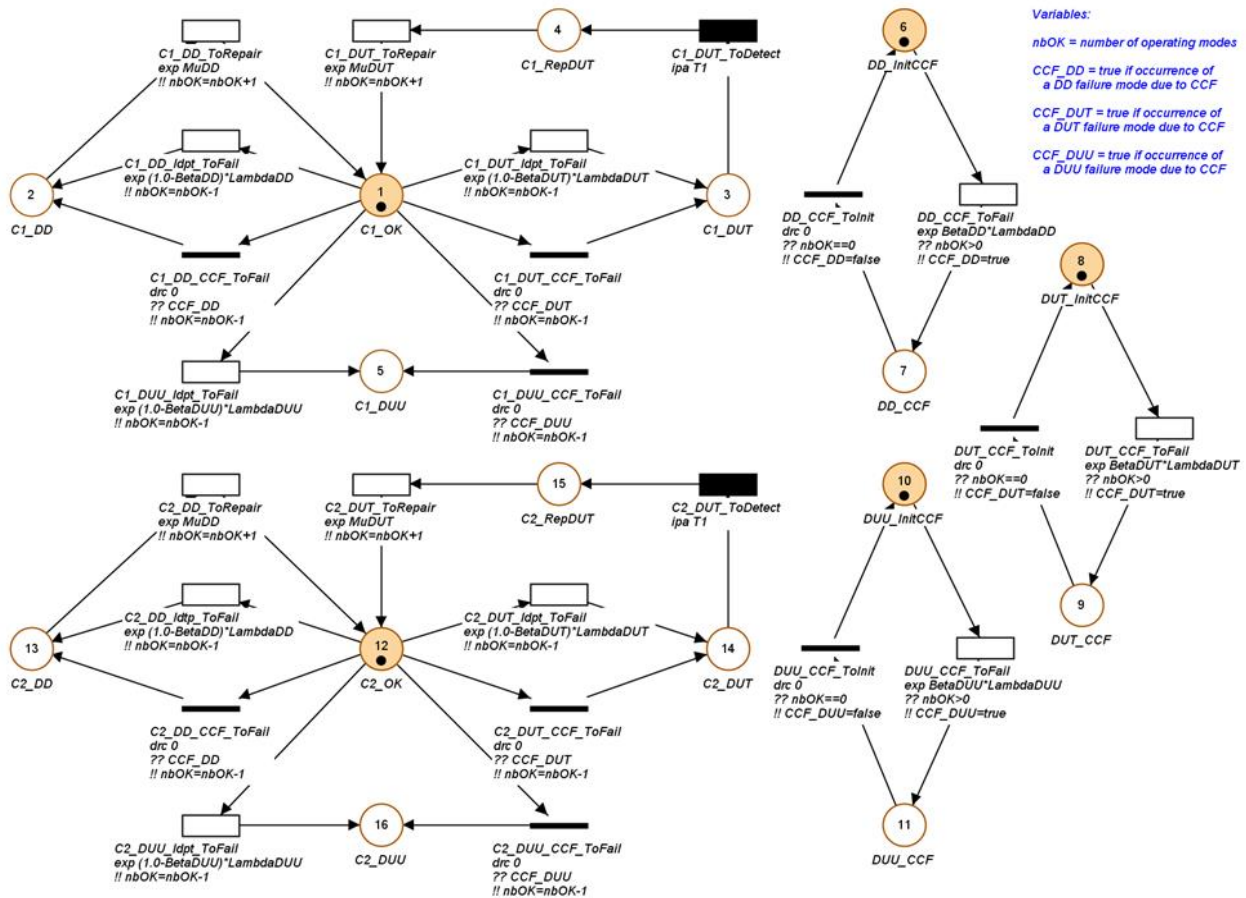

**Figure 3.** Réseaux de Petri stochastiques à prédicats avec GRIF/Petri, cas d'étude avec une architecture *1oo2*

En utilisant GRIF/Petri, la première ligne sous une place ou sous une transition est la description de l'élément. Par exemple, lorsqu'une place avec la description « Cx_y » (pour les Places 1 à 5 et 12 à 16) contient un jeton, cela signifie que le Canal *x* est dans le mode opérant, de défaillance, ou de réparation *y* (suivant les notations données dans la Section II.2). (Ici, il est considéré que chaque canal ne peut pas se trouver dans les modes multiples au même instant, tout comme pour le modèle Markovien, cf. Section III.2.) Une place qui contient au moins un jeton est représentée en couleur (pour ce cas d'étude, une place ne peut contenir que *0* ou *1* jeton). Sous une transition, l'expression commençant par « ?? » est une « garde », et l'expression commençant par « !! » est une « affectation ». Par exemple, la transition de la Place 6 à la Place 7 est franchissable uniquement si la variable « nbOK » est plus grande que *0* (i.e. au moins un canal est dans un mode opérant) et, lorsqu'elle est franchie, la variable « CCF_DD » prend la valeur « true » (i.e. occurrence d'un mode de défaillance « DD » due à une cause commune). Puisque « CCF_DD » est « true », la transition de la Place 1 à la Place 2 (respectivement de la Place 12 à la Place 13) est donc franchie si la Place 1 (respectivement, la Place 12) contient un jeton. Si cette transition est franchie, la variable « nbOK » est réduite de *1* (i.e. le nombre de canaux dans un mode opérant est réduit de *1*). Si « nbOK » devient égal à *0* (i.e. aucun canal n'est dans un mode opérant), alors la transition de la Place 7 à la Place 6 rétablit la variable « CCF_DD » à « false ». En utilisant GRIF/Petri, une transition avec un délai distribué exponentiellement est représentée en blanc, une transition sans délai est représentée par un fin rectangle noir, et tout autre transition déterministe est représentée par un large rectangle noir. Pour les délais de transition, trois types de lois sont ici utilisés :

- « exp », qui correspond à une distribution exponentielle, définie par un taux (e.g. taux de défaillance, taux de réparation) ;
- « drc », qui correspond à une distribution de Dirac (i.e. délai de transition déterministe), définie par un délai constant (ici égal à *0* afin d'avoir des transitions instantanées dès que les « gardes » sont vérifiées) ;
- « ipa », qui signifie « *in advance schedule times* », définie par une période où la transition devient ponctuellement franchissable (et qui est utilisée avec les périodes de tests de révision afin de détecter les modes de défaillance « DUT » à des instants déterministes). À noter que contrairement aux autres lois, celle-ci est basée sur le temps calendaire, et non sur la durée pendant laquelle la transition est restée franchissable.

D'autres lois sont également disponibles dans GRIF/Petri, incluant les lois : uniforme, triangulaire, Weibull, normale, et log-normale. De plus, et parce que ces réseaux de Petri permettent des modélisations comportementales combinées à des analyses basées sur des simulations, ils fournissent un moyen de modéliser d'avantage de caractéristiques (toutes ?) telles que des tests asynchrones, des tests et/ou des réparations imparfaites, des défaillances générées par des tests et/ou des réparations, des indisponibilités et/ou des modes dégradés durant des tests, etc.

### III.4. Équations approchées avec OrbitSIL

Les équations approchées sont basées sur des développements limités (i.e. séries de Taylor), (souvent du premier ordre), de fonctions exponentielles afin d'approcher les équations de fiabilité par des formules qui ne contiennent que des opérateurs « basiques » (e.g. sans intégrales), et d'être ainsi simplement calculable. Comme les développements limités du premier ordre surestiment les fonctions exponentielles de défiabilité, ces équations approchées sont souvent considérées comme « conservatives » (e.g. surestimation de la *PFD*$_{avg}$). Cependant, il convient de préciser que ces équations sont développées avec des hypothèses sous-jacentes spécifiques et que les appliquer en dehors de leurs limites peut conduire à des résultats (dangereusement) erronés (e.g. lorsque les tests de révision ne sont pas efficace à *100*%).

Les formules ci-dessous ont été développées par DNV et sont différentes de celles données dans la norme CEI 61508 (Partie 6) (CEI, 2010). Ces formules ont été proposées par L.F. Oliveira *et al.*, notamment avec l'inclusion de tests imparfaits (Oliveira, 2008) et de défaillances de causes communes (Oliveira & Abramovitch, 2010), afin d'être plus complètes que celles de la norme CEI 61508. Cependant, ces formules sont aussi plus « compliquées » et le logiciel (non-commercialisé) OrbitSIL a donc été développé pour les exploiter simplement. Parce qu'elles sont basées sur des développements limités du premier ordre, les équations approchées suivantes, ici utilisées pour calculer les *PFD*$_{avg}$ du cas d'étude, ne sont considérées comme « valides » uniquement si $\lambda_{DUT} \times T_1 < 0.1$ et $\lambda_{DUU} \times T_0 < 0.1$ (cf. notations et hypothèses données dans la Section II):

$$
\begin{aligned}
PFD_{avg} =& \binom{N}{N-M+1} \times \left( (1-\beta_{DD}) \times \frac{\lambda_{DD}}{\mu_{DD}} \right)^{N-M+1} \\
&+ \binom{N}{N-M+1} \times \left( (1-\beta_{DUT}) \times \lambda_{DUT} \right)^{N-M+1} \times T_1^{N-M} \times \left( \frac{T_1}{N-M+2} + \frac{1}{\mu_{DUT}} \right) + \binom{N}{N-M+1} \times \left( (1-\beta_{DUU}) \times \lambda_{DUU} \right)^{N-M+1} \times T_0^{N-M} \times \left( \frac{T_0}{N-M+2} \right) \\
&+ \sum_{i=1}^{N-M} \binom{N-i}{N-M+1-i} \times \left( f(N-M+1-i, \beta_{DD}) \times \frac{\lambda_{DD}}{\mu_{DD}} \right)^{N-M+1-i} \times \binom{N}{i} \times (f(i, \beta_{DUT}) \times \lambda_{DUT})^i \times T_1^{i-1} \times \left( \frac{T_1}{i+1} + \frac{1}{\mu_{DUT}} \right) \\
&+ \sum_{i=1}^{N-M} \binom{N-i}{N-M+1-i} \times \left( f(N-M+1-i, \beta_{DD}) \times \frac{\lambda_{DD}}{\mu_{DD}} \right)^{N-M+1-i} \times \binom{N}{i} \times (f(i, \beta_{DUU}) \times \lambda_{DUU})^i \times T_0^{i-1} \times \left( \frac{T_0}{i+1} \right) \\
&+ \sum_{i=1}^{N-M} \binom{N-i}{N-M+1-i} \times \left( (1-\beta_{DUT}) \times \lambda_{DUT} \right)^{N-M+1-i} \times T_1^{N-M-i} \times \left( \frac{T_1}{N-M+2-i} + \frac{1}{\mu_{DUT}} \right) \times \binom{N}{i} \times \left( (1-\beta_{DUU}) \times \lambda_{DUU} \right)^i \times T_0^{i-1} \times \left( \frac{T_0}{i+1} \right) \\
&+ \sum_{i=1}^{N-M-1} \sum_{j=1}^{N-M-i} \binom{N-i-j}{N-M+1-i-j} \times \left( f(N-M+1-i-j, \beta_{DD}) \times \frac{\lambda_{DD}}{\mu_{DD}} \right)^{N-M+1-i-j} \\
&\times \binom{N-i}{j} \times \left( (1-\beta_{DUT}) \times \lambda_{DUT} \right)^j \times T_1^{j-1} \times \left( \frac{T_1}{j+1} + \frac{1}{\mu_{DUT}} \right) \times \binom{N}{i} \times \left( (1-\beta_{DUU}) \times \lambda_{DUU} \right)^i \times T_0^{i-1} \times \left( \frac{T_0}{i+1} \right) \\
&+ \beta_{DD} \times \frac{\lambda_{DD}}{\mu_{DD}} + \beta_{DUT} \times \lambda_{DUT} \times \left( \frac{T_1}{2} + \frac{1}{\mu_{DUT}} \right) + \beta_{DUU} \times \lambda_{DUU} \times \left( \frac{T_0}{2} \right)
\end{aligned}
\quad \{1\}
$$

avec la fonction *f(x,b)* définie par :

$$f(x,b) = \begin{cases} 1, & \text{if } x = 1 \\ 1 - b, & \text{if } x \neq 1 \end{cases}$$

et la fonction combinatoire suivante:

$$\binom{n}{k} = \frac{n!}{k! \times (n-k)!}$$

En accord avec les hypothèses de la Section II.1, l'équation approchée {1} ne prend pas en compte des temps de réparation pour les défaillances dangereuses uniquement détectées par des demandes réelles, mais peuvent néanmoins être facilement étendues afin d'y inclure ce paramètre. D'autres caractéristiques, comme des taux de défaillance et de réparation hétérogènes, des tests asynchrones, et des défaillances générées par des tests et/ou des réparations, sont quant à elles plus difficile à intégrer dans de telles équations. Sous d'autres hypothèses, et notamment lorsque les temps de réparation n'ont pas besoin d'être modélisés (e.g. lorsque des mesures compensatoires sont mises en place afin de maintenir la fonction de sécurité pendant les réparations), d'autres équations « approchées » et « exactes » ont été développées, qui peuvent notamment inclure des tests non-périodiques (Brissaud *et al.*, 2010).

# IV. Discussions

### IV.1. Résultats

Les résultats obtenus par l'application de chaque méthode, pour les différentes configurations du cas d'étude présenté dans la Section II, sont reportés dans la Table 2. Pour les cas *v* et *vi*, les graphes de Markov auraient contenus 35 états, ce qui a découragé les auteurs d'appliquer les modèles Markovien à ces deux configurations. Pour les réseaux de Petri, $10^8$ simulations ont été utilisées pour chaque cas, permettant ainsi d'obtenir des intervalles de confiance à 90% allant de $\pm 0.03\%$ (pour le cas *ii*) à $\pm 0.5\%$ (pour le cas *iii*). Le temps de simulation a alors varié de 40 minutes (pour les cas *i* et *ii*) à 75 minutes (pour les cas *v* et *vi*), avec un processeur de 2,67 GHz et 4,00 GB de RAM.

**Table 2.** Résultats : $PFD_{avg}$ pour les cas d'étude *i* à *vi*

|  | cas *i* | cas *ii* | cas *iii* | cas *iv* | cas *v* | cas *vi* |
|---|---|---|---|---|---|---|
| **Arbres de défaillance** | $7.43 \times 10^{-3}$ | $1.27 \times 10^{-1}$ | $4.31 \times 10^{-4}$ | $2.93 \times 10^{-2}$ | $5.48 \times 10^{-4}$ | $5.59 \times 10^{-2}$ |
| **Modèles Markoviens** | $7.41 \times 10^{-3}$ | $1.24 \times 10^{-1}$ | $4.29 \times 10^{-4}$ | $2.83 \times 10^{-2}$ | - | - |
| **Réseaux de Petri** | $7.41 \times 10^{-3}$ | $1.24 \times 10^{-1}$ | $4.30 \times 10^{-4}$ | $2.83 \times 10^{-2}$ | $5.47 \times 10^{-4}$ | $5.43 \times 10^{-2}$ |
| **Équations approchées** | $7.46 \times 10^{-3}$ | $1.38 \times 10^{-1}$ | $4.31 \times 10^{-4}$ | $3.25 \times 10^{-2}$ | $5.49 \times 10^{-4}$ | $6.98 \times 10^{-2}$ |

Les graphes de Markov multi-phase et les réseaux de Petri stochastiques à prédicats produisent les mêmes résultats (la différence pour le cas *iii* n'est pas significative compte tenu des intervalles de confiance) car les mêmes hypothèses ont été faites. Les résultats obtenus par les arbres de défaillance sont légèrement plus grands (avec un écart maximum de *3,5%* pour le cas *iv*) de par l'indépendance des évènements de base, qui autorise des modes de défaillance multiples du même canal au même instant (ces combinaisons ayant volontairement été exclues dans les modèles Markoviens et les réseaux de Petri utilisés). Enfin, les équations approchées sont celles qui surestiment le plus les résultats (avec un écart maximum de *28,5%* pour le cas *vi*, comparé aux réseaux de Petri) de par les approximations mathématiques qui ont été faites pour disposer de formules « simples » (à noter que pour les cas *ii*, *iv*, et *vi*, $\lambda_{DUU} \times T_0 = 0.213$, ce qui est plus grand que *0.1*, et ce qui ne satisfait donc pas aux « conditions de validité » définies dans la Section III.4).

Pour résumer, différentes méthodes fondamentalement différentes (à la fois d'un point de vue pratique et mathématique) sont capables de produire des résultats similaires pour l'estimation de la probabilité moyenne de défaillance dangereuse en cas de sollicitation ($PFD_{avg}$), en tenant compte des caractéristiques requises par la norme CEI 61508. Pour cela, il est néanmoins nécessaire de maîtriser la méthode utilisée (ce qui implique notamment de connaître les hypothèses intrinsèques et les limites de modélisation de cette méthode), et d'utiliser un outil logiciel adapté et performant (ce qui implique notamment qu'il exploite pleinement et correctement les théories de la fiabilité). Le choix d'une méthode doit alors se faire sans a priori, sur la base d'un accord entre les efforts de modélisation, les objectifs, et les propriétés du système. Afin d'assister l'analyste dans cette tâche, une discussion sur les avantages et limites de chacune de ces méthodes est proposé ci-après.

### IV.2. Synthèse comparative des méthodes

La Table 3 fournit quelques critères (qui ne prétendent pas à être exhaustifs) pour comparer les arbres de défaillance, les graphes de Markov multi-phase, les réseaux de Petri stochastiques à prédicats, et les équations approchées, tels que présentés dans la Section III (à noter que ces méthodes possèdent aussi d'autres extensions qui ne sont pas discutées ici). Cette « évaluation comparative » est basée sur l'expérience de ses auteurs, et reste bien sûr ouverte à la discussion.

# V. Conclusions

- Les équations approchées fournissent un moyen simple et rapide d'évaluer de nombreux systèmes simples/basiques. Cependant, cette approche est aussi, par nature, la moins flexible et il peut être hasardeux de l'utiliser en dehors des hypothèses sous-jacentes. L'inconvénient de ce type d'approche est ainsi le recours potentiel de certains utilisateurs à des formules « toutes faites » qui paraissent donc « attractives », sans prendre en considération les précautions requises (notamment dans le cas d'hypothèses non « conservatives »).
- Les arbres de défaillance possèdent la plupart des avantages pour les ingénieurs : relativement simples à appliquer et à lire, basés sur des méthodes performantes d'analyses, et qui permettent l'inclusion de nombreuses caractéristiques (dont des caractéristiques avancées de test et de réparation) – à condition bien sûr d'utiliser d'un outil logiciel performant. La seule limite concerne les caractéristiques « dynamiques » (e.g. séquences d'évènements) qui ne peuvent pas être modélisées à cause d'hypothèses intrinsèques aux arbres de défaillance (e.g. indépendance entre les évènements de base).

- Les réseaux de Petri stochastiques à prédicats fournissent l'approche la plus flexible pour modéliser des systèmes complexes, notamment ceux qui possèdent des propriétés particulières telles que des caractéristiques dynamiques. Cette flexibilité est permise par l'approche d'analyse qui est basée sur des simulations. La contrepartie à cela est qu'un temps parfois important est requis pour effectuer les analyses. Les réseaux de Petri devraient donc être recommandés pour les cas qui ne peuvent pas être traités convenablement par des arbres de défaillance.
- Les graphes de Markov multi-phase n'ont pas de réels avantages comparés aux autres méthodes, et elle est probablement l'approche qui a le ratio coût-bénéfice le plus faible en termes d'efforts de modélisation versus capacités de modélisation.

**Table 3.** Comparaison* des méthodes

|  | **Arbres de défaillance** | **Modèles Markoviens** | **Réseaux de Petri** | **Équations approchées** |
|---|---|---|---|---|
| **Modèle** | | | | |
| taille du modèle (pour de grands systèmes) | linéairement dépendant | exponentiellement dépendant | linéairement dépendant | fixe |
| temps de modélisation (pour de nombreux systèmes) | relativement rapide | assez long | assez long mais peut être réduit par l'utilisation de prototypes | rapide (une seule formule à répéter) |
| assujetti aux erreurs de modélisation | modèles « simples » / nombre limité d'erreurs possibles | modèles « riches » / nombreuses erreurs possibles | modèles « riches » / nombreuses erreurs possibles | modèles basiques / peu d'erreurs possibles |
| lisibilité du modèle | simple pour des ingénieurs | réservée aux spécialistes | réservée aux spécialistes | « boite noire » difficile à lire |
| **Analyses** | | | | |
| méthode d'analyse | approche Booléenne « exacte » | approche mathématique « exacte » | simulations de Monte Carlo | résultats approchés |
| temps d'analyse / analyses d'incertitudes | très rapide avec de bons algorithmes | très rapide avec de bons algorithmes | assez long à cause des simulations | très rapide |
| adapté aux analyses de facteurs d'importance | plusieurs facteurs d'importance dédiés | à adapter | à adapter | non adapté |
| **Flexibilité** | | | | |
| taux de défaillance et/ou de réparation hétérogènes | bonnes capacités | capable mais accroisse drastiquement le modèle | bonnes capacités | non adapté |
| caractéristiques avancées de test et/ou de réparation | bonnes capacités avec un logiciel approprié | capable mais accroisse drastiquement le modèle | bonnes capacités | non adapté |
| liberté des distributions (e.g. Weibull, normale, uniforme) | bonnes capacités avec un logiciel approprié | non adapté | bonnes capacités | non adapté |
| caractéristiques dynamiques (e.g. séquences) | non adapté | capable pour des cas basiques | bonnes capacités | non adapté |

*Dans cette table, la « souhaitabilité » de chaque méthode au regard des critères proposés est classée de « vert » (i.e. la plus souhaitable et/ou la moins restrictive) à « rouge » (i.e. la moins souhaitable et/ou la plus restrictive). Cette classification est basée sur l'expérience de ses auteurs, et reste bien sûr ouverte à la discussion.

# VI. **Références**